\documentclass{article}

\usepackage{amssymb,amsmath,theorem,euscript}

\input{epsf}

\newcounter{sec}

\newcounter{punct}[sec]

\def\punct{\refstepcounter{punct}{\arabic{sec}.\arabic{punct}.  }\boldmath}

\def\COUNTERS{\addtocounter{sec}{1}
              \setcounter{punct}{0}
          \setcounter{equation}{0}
          \setcounter{theorem}{0}
                  }

\newtheorem{theorem}{Theorem}[sec]
\newtheorem{proposition}[theorem]{Proposition}

\newtheorem{lemma}[theorem]{Lemma}

 \def\ov{\overline}
\def\wt{\widetilde}

\renewcommand{\Im}{\mathop {\mathrm {Im}}\nolimits}

 \newcommand{\indef}{\mathop {\mathrm {indef}}\nolimits}
 \newcommand{\dom}{\mathop {\mathrm {dom}}\nolimits}

\begin{document}

\def\OO{\mathrm{O}}
\def\GLO{\mathrm{GLO}}
\def\Coll{\mathrm{Coll}}
\def\kappa{\varkappa}
\def\Mat{\mathrm{Mat}}
\def\U{\mathrm U}
\def\GLL{\overline{\mathrm {GL}}}
\def\Spp{\overline{\mathrm {Sp}}}
\def\SL{\mathrm{SL}}
\def\PGL{\mathrm{PGL}}

\def\R{\mathbb{R}}
\def\C{\mathbb{C}}
\def\V{\mathbb{V}}

\def\la{\langle}
\def\ra{\rangle}

 \def\cA{\mathcal A}
\def\cB{\mathcal B}
\def\cC{\mathcal C}
\def\cD{\mathcal D}
\def\cE{\mathcal E}
\def\cF{\mathcal F}
\def\cG{\mathcal G}
\def\cH{\mathcal H}
\def\cJ{\mathcal J}
\def\cI{\mathcal I}
\def\cK{\mathcal K}
 \def\cL{\mathcal L}
\def\cM{\mathcal M}
\def\cN{\mathcal N}
 \def\cO{\mathcal O}
\def\cP{\mathcal P}
\def\cQ{\mathcal Q}
\def\cR{\mathcal R}
\def\cS{\mathcal S}
\def\cT{\mathcal T}
\def\cU{\mathcal U}
\def\cV{\mathcal V}
 \def\cW{\mathcal W}
\def\cX{\mathcal X}
 \def\cY{\mathcal Y}
 \def\cZ{\mathcal Z}
 
 \def\cGL{\mathcal{GL}}
 
 \def\EXP{\mathcal{EXP}}
 
 \def\frA{\mathfrak A}
 \def\frB{\mathfrak B}
\def\frC{\mathfrak C}
\def\frD{\mathfrak D}
\def\frE{\mathfrak E}
\def\frF{\mathfrak F}
\def\frG{\mathfrak G}
\def\frH{\mathfrak H}
\def\frI{\mathfrak I}
 \def\frJ{\mathfrak J}
 \def\frK{\mathfrak K}
 \def\frL{\mathfrak L}
\def\frM{\mathfrak M}
 \def\frN{\mathfrak N} \def\frO{\mathfrak O} \def\frP{\mathfrak P} \def\frQ{\mathfrak Q} \def\frR{\mathfrak R}
 \def\frS{\mathfrak S} \def\frT{\mathfrak T} \def\frU{\mathfrak U} \def\frV{\mathfrak V} \def\frW{\mathfrak W}
 \def\frX{\mathfrak X} \def\frY{\mathfrak Y} \def\frZ{\mathfrak Z} \def\fra{\mathfrak a} \def\frb{\mathfrak b}
 \def\frc{\mathfrak c} \def\frd{\mathfrak d} \def\fre{\mathfrak e} \def\frf{\mathfrak f} \def\frg{\mathfrak g}
 \def\frh{\mathfrak h} \def\fri{\mathfrak i} \def\frj{\mathfrak j} \def\frk{\mathfrak k} \def\frl{\mathfrak l}
 \def\frm{\mathfrak m} \def\frn{\mathfrak n} \def\fro{\mathfrak o} \def\frp{\mathfrak p} \def\frq{\mathfrak q}
 \def\frr{\mathfrak r} \def\frs{\mathfrak s} \def\frt{\mathfrak t} \def\fru{\mathfrak u} \def\frv{\mathfrak v}
 \def\frw{\mathfrak w} \def\frx{\mathfrak x} \def\fry{\mathfrak y} \def\frz{\mathfrak z} \def\frsp{\mathfrak{sp}}
 \def\bfa{\mathbf a} \def\bfb{\mathbf b} \def\bfc{\mathbf c} \def\bfd{\mathbf d} \def\bfe{\mathbf e} \def\bff{\mathbf f}
 \def\bfg{\mathbf g} \def\bfh{\mathbf h} \def\bfi{\mathbf i} \def\bfj{\mathbf j} \def\bfk{\mathbf k} \def\bfl{\mathbf l}
 \def\bfm{\mathbf m} \def\bfn{\mathbf n} \def\bfo{\mathbf o} \def\bfp{\mathbf p} \def\bfq{\mathbf q} \def\bfr{\mathbf r}
 \def\bfs{\mathbf s} \def\bft{\mathbf t} \def\bfu{\mathbf u} \def\bfv{\mathbf v} \def\bfw{\mathbf w} \def\bfx{\mathbf x}
 \def\bfy{\mathbf y} \def\bfz{\mathbf z} \def\bfA{\mathbf A} \def\bfB{\mathbf B} \def\bfC{\mathbf C} \def\bfD{\mathbf D}
 \def\bfE{\mathbf E} \def\bfF{\mathbf F} \def\bfG{\mathbf G} \def\bfH{\mathbf H} \def\bfI{\mathbf I} \def\bfJ{\mathbf J}
 \def\bfK{\mathbf K} \def\bfL{\mathbf L} \def\bfM{\mathbf M} \def\bfN{\mathbf N} \def\bfO{\mathbf O} \def\bfP{\mathbf P}
 \def\bfQ{\mathbf Q} \def\bfR{\mathbf R} \def\bfS{\mathbf S} \def\bfT{\mathbf T} \def\bfU{\mathbf U} \def\bfV{\mathbf V}
 \def\bfW{\mathbf W} \def\bfX{\mathbf X} \def\bfY{\mathbf Y} \def\bfZ{\mathbf Z} \def\bfw{\mathbf w}
\def\B{\mathrm{B}}
 \def\R {{\mathbb R }} \def\C {{\mathbb C }} \def\Z{{\mathbb Z}} \def\H{{\mathbb H}}
  \def\K{{\mathbb K}}
   \def\k{{\Bbbk}}
 \def\N{{\mathbb N}} \def\Q{{\mathbb Q}} \def\A{{\mathbb A}} \def\T{\mathbb T} 
 \def\G{\mathbb G}
 \def\bbA{\mathbb A} \def\bbB{\mathbb B} \def\bbD{\mathbb D} \def\bbE{\mathbb E} \def\bbF{\mathbb F} \def\bbG{\mathbb G}
 \def\bbI{\mathbb I} \def\bbJ{\mathbb J} \def\bbL{\mathbb L} \def\bbM{\mathbb M} \def\bbN{\mathbb N} \def\bbO{\mathbb O}
 \def\bbP{\mathbb P} \def\bbQ{\mathbb Q} \def\bbS{\mathbb S} \def\bbT{\mathbb T} \def\bbU{\mathbb U} \def\bbV{\mathbb V}
 \def\bbW{\mathbb W} \def\bbX{\mathbb X} \def\bbY{\mathbb Y} \def\kappa{\varkappa} \def\epsilon{\varepsilon}
 \def\phi{\varphi} \def\le{\leqslant} \def\ge{\geqslant}

\def\P{\mathrm P}
\def\T{\mathrm T}

\def\GL{\mathrm {GL}}
\def\bGL{\mathbf {GL}}
\def\GLB{\mathrm {GLB}}

\def\bGr{\mathbf {Gr}}
\def\Gr{\mathrm {Gr}}
\def\Sp{\mathrm {Sp}}
\def\bFl{\mathbf {Fl}}

\def\1{\mathbf {1}}
\def\0{\mathbf {0}}

\def\rra{\rightrightarrows}

 \newcommand{\Dim}{\mathop {\mathrm {Dim}}\nolimits}
  \newcommand{\codim}{\mathop {\mathrm {codim}}\nolimits}
   \newcommand{\im}{\mathop {\mathrm {im}}\nolimits}
\newcommand{\ind}{\mathop {\mathrm {ind}}\nolimits}
\newcommand{\graph}{\mathop {\mathrm {graph}}\nolimits}
\newcommand{\Hinge}{\mathop {\mathrm {Hinge}}\nolimits}
\newcommand{\hinge}{\mathop {\mathrm {hinge}}\nolimits}

\def\F{\bbF}

\def\lambdA{{\boldsymbol{\lambda}}}
\def\alphA{{\boldsymbol{\alpha}}}
\def\betA{{\boldsymbol{\beta}}}
\def\gammA{{\boldsymbol{\gamma}}}
\def\deltA{{\boldsymbol{\delta}}}
\def\mU{{\boldsymbol{\mu}}}
\def\nU{{\boldsymbol{\nu}}}
\def\epsiloN{{\boldsymbol{\varepsilon}}}
\def\phI{{\boldsymbol{\phi}}}
\def\psI{{\boldsymbol{\psi}}}
\def\kappA{{\boldsymbol{\kappa}}}

\def\sm{\smallskip}
\def\nw{\nwarrow}
\def\se{\searrow}

\begin{center}
	\Large \bf
	On double cosets of groups $\GL(n)$ with respect to subgroups of block
	strictly triangular matrices
	
	\bigskip 
	
	\large\sc
	Yury A. Neretin\footnote{Supported by the grant of  FWF (Austrian Science Fund), P31591.}
\end{center}

\bigskip

{\small We parametrize the space of double cosets of the group $\mathrm{GL}(n,\Bbbk)$ with respect
	to two subgroups $\mathrm{T}_-$, $\mathrm{T}_+$ of block 	strictly triangular matrices.
In Addendum, we consider the quasi-regular representation of $\GL(n,\mathbb{C})$ in  $L^2$ on 
$\mathrm{T}_-\setminus\GL(n,\mathbb{C})$, observe that it admits an additional group of symmetries,
find the joint spectrum, and observe that it is multiplicity free.}

\section{The statement}

\COUNTERS 

{\bf \punct Double cosets}
Let $G$ be a group, $K$, $L$ subgroups. A {\it double coset} of $G$
with respect to $K$, $L$ is a set of the type $K\cdot g \cdot L$,
i.e., the set of all elements of $G$ that can be represented in the form
$kgl$, where $g$ is fixed, $k$ ranges in $K$, $l$ ranges in $L$.
We denote the set of all double cosets by $K\setminus G/L$.

 A description
of this set is equivalent to a description of orbits of $K$ on the homogeneous space
$G/L$, and to a description of orbits of $L$ on the homogeneous space
$K\setminus G$. If $G$ is finite, then a description of double cosets
is equivalent to a description of intertwining operators between quasi-regular
representations of $G$ in $\ell^2(G/K)$ and $\ell^2(G/L)$
(see, e.g. \cite{Kir}, Sect. 13.1). For Lie groups and locally compact groups
picture is more complicate, in any case understanding double cosets
  seems necessary for understanding  analysis on the corresponding homogeneous spaces.

 In any case a problem of description of double cosets arises quite often, but not quite often it admits a 
tame solution.

\sm

{\bf\punct The problem.}
Let $\k$ be a field, $V$ be a finite dimensional linear space over $\k$.
Denote by $\GL[V]$ be the group of all invertible linear operator in $V$.
We also use the notation
 $\GL(n)=\GL(n,\k)$ for $\GL[\k^n]$.  Split $V$ into a direct sum
 $$
 V=V_1\oplus \dots\oplus V_p, \qquad \dim V_j=\alpha_j.
 $$
 Denote 
 $$\T_+[V_1, \dots, V_p]=\T_+(\alpha_1,\dots,\alpha_p)=\T_+(|\alphA|)
 $$
 the group of all block strictly upper upper triangular matrices of the size
 $(\alpha_1+\dots+\alpha_p)\times (\alpha_1+\dots+\alpha_p)$, i.e.,
 matrices of the form
 \begin{equation}
 \begin{pmatrix}
 1_{\alpha_1}&*&\dots&*\\
 0&1_{\alpha_2}&\dots &*\\
 \vdots &\vdots &\ddots &\vdots\\
 0&0&\dots&1_{\alpha_p}
 \end{pmatrix},
 \end{equation}
 where $1_m$ denotes the unit matrix of size $m$.
 By $\P_+[\dots]=\P_+(\dots)$ we denote the group of all block triangular
 matrices, i.e., we allow arbitrary invertible matrices 
 on places of units. Clearly, $\T_+$ is normal in $\P_+$,
 $$
 \P_+(\alpha_1,\dots,\alpha_p)/\T_+(\alpha_1,\dots,\alpha_p)
 \simeq \prod_{i=1}^p \GL(\alpha_i).
 $$
 By $\T_-[\dots]=\T_-(\dots)$ and $\P_-[\dots]=\P_-(\dots)$
 we denote the corresponding groups of  lower triangular matrices.

In this paper we describe double coset spaces
\begin{equation}
\T_-(\beta_1,\dots,\beta_q)\setminus \GL(n)/\T_+(\alpha_1,\dots,\alpha_p).
\label{eq:problem}
\end{equation}

 {\bf\punct The statement.} Recall some definitions. Let $X$, $Y$ be linear spaces
 over $\k$. A {\it linear relation} (see, e.g., \cite{Ner-book}, Sect. 2.5) $L:X\rra Y$ is a linear
 subspace in $X\oplus Y$. A graph of a linear map $X\to Y$ is a linear relation
 but not vice versa. For a linear relation we define: 
 
 \sm 
 
 1) the {\it kernel} $\ker L$ is the intersection $L\cap X$;
 
 \sm 
 
 2) the {\it image} $\im L$ is the projection of $L$ to $Y$ along $X$;
 
 \sm 
 
 3) the {\it domain} $\dom L$ is the projection of $L$ to $X$ along $Y$;
 
 \sm 
 
 4) the {\it indefiniteness} $\indef L$ is the intersection $L\cap Y$;
 
 \sm 
 

If $L$ is a graph of a linear operator $A:X\to Y$, then
kernel, image, and rank are the usual kernel, image, and rank.
In this case also $\dom L=X$, $\indef L=0$.

Any linear relation $L$ determines a canonical invertible operator
$$
\Theta(L):\dom L/\ker L\to \im L/\indef L.
$$
Moreover a linear relation $L:X\rra Y$ is determined by subspaces
$\ker L\subset \dom L$ in $X$, subspaces $\indef L\subset \im L$ in
$Y$ and the operator $\Theta(L)$.

\sm 

Now consider a linear space $V\simeq \K^n$ and take two decompositions
\begin{align*}
V:=\k^{\alpha_1}\oplus\dots \oplus \k^{\alpha_p}=:V_1\oplus\dots\oplus V_p;
\\
W:=\k^{\beta_1}\oplus\dots \oplus \k^{\beta_p}=:W_1\oplus\dots\oplus W_q.
\end{align*}
Consider the double cosets (\ref{eq:problem}). For each element $A\in \GL[V]$
we assign a canonical collection of linear relations 
$$
\chi_{ij}(A):V_i\rra W_j
$$
defined in the following way. We say that $(\xi, \eta)\in \chi_{ij}(A)$
if there exist $x_1\in V_1$, \dots, $x_{i-1}\in V_{i-1}$
and $y_{j+1}\in W_{j+1}$, \dots, $y_q\in W_q$ such that
\begin{equation*}
\begin{pmatrix}
0\\ \vdots\\0\\ \eta\\y_{j+1}\\ \vdots\\ y_p
\end{pmatrix}=
A
\begin{pmatrix}
x_1\\ \vdots \\ x_{i-1}\\ \xi\\0\\ \vdots\\0
\end{pmatrix}
\end{equation*}

\begin{proposition}
	\label{pr:}
	{\rm a)}
	The relations $\chi_{ij}$ depend only on the double coset containing $A$;
	
	\sm 
	
	{\rm b)} The relations $\chi_{ij}=\chi_{ij}(A)$ satisfy
	the conditions:
	\begin{equation}
	\ker \chi_{ij}=\dom \chi_{i(j+1)},\qquad
	\im\chi_{ij}=\indef \chi_{(i+1)j},
	\label{eq:chiij}
	\end{equation}
	and
	\begin{align}
	\indef \chi_{1j}=0,\qquad \im \chi_{pj}=W_j;
	\label{eq:kraj1}
	\\
	\ker \chi_{iq}=0, \qquad \dom \chi_{i1}=V_i.
	\label{eq:kraj2}
	\end{align}
\end{proposition}

We say that a collection of linear relations $\xi_{ij}:V_i\rra W_j$ is a {\it bi-hinge%
\footnote{Cf. similar objects in \cite{Ner-universal}.}}
if it satisfies conditions
(\ref{eq:chiij})--\eqref{eq:kraj2}. We denote the set of all hinges
by 
$$\Hinge=\Hinge[V_1,\dots,V_p;W_1,\dots,W_q].$$

\begin{theorem}
	\label{th:}
	The map $A\mapsto \chi_{ij}(A)$ determines a one-to-one correspondence between
	the double coset space \eqref{eq:problem} and the set of bi-hinges.
\end{theorem}

{\bf \punct Some known cases of classification of double cosets.}
Discuss shortly such cases related to classical (also, semisimple, reductive) groups.
The are 3 big important series of solvable problems with tame solutions.

\sm 

---  Let $G$ be a real semisimple group, $H$, $L$ are symmetric subgroups%
\footnote{Recall that a subgroup $H\subset G$ is symmetric if it is
the set of fixed points of some involution $\sigma:G\to G$ 
(i.e., $\sigma(g_1)\sigma(g_2)=\sigma(g_2g_1)$, $\sigma(\sigma(g))=g$.}. In particular, this problem
includes the Jordan normal form (more generally, description of conjugacy classes in all semisimple Lie groups%
\footnote{Namely, we set $G=Q\times Q$, $H=L=\mathrm{diag}\,Q$ is the diagonal 
subgroup.} $Q$),  reduction of pairs of nondegenerate quadratic (or symplectic) forms, 
canonical forms of pairs of subspaces  in a Euclidean space, etc. A formal reference
to a 'general case' is \cite{Mat-symm}.

\sm

--- We consider $H\setminus G/P$, where $G$ is semisimple group, $H$ is a symmetric subgroup, $P$ is a parabolic subgroup
(a block triangular subgroup). A formal reference to a 'general case' is \cite{Mat-parab}.

\sm 

--- For $p$-adic groups the most important case is related to the Iwahori
 subgroups, see \cite{Iwa}.
 
 \sm 
 
 There is a cloud of minor variations of these series (we can slightly enlarge $G$,
  or slightly reduce subgroups).
  
  \sm

  Next, there are different ways to assign spectral data to several matrices
  (this also an be regarded as a classification of double cosets):
  a spectral curve with a bundle, see \cite{Tyu}, \cite{CT}, \cite{Hit},
  or a spectral surface with a sheaf, see \cite{Beu}.
  
  On the other hand, for infinite-dimensional groups quite often 
  a double coset space $K\setminus G/K$ has a structure of a semigroup.
  There arise questions about spectral data visualizing such multiplications.
  This also leads to objects of algebraic-geometric nature
  as spaces of holomorphic maps of Riemann sphere to Grassmannians
  (see \cite{Ner-book}, Sect.X.3) or rational maps of Grassmannians to Grassmannians
  (see \cite{Ner-char}). 
  
  Our case arose as  a byproduct of a construction
  of the latter type in \cite{Ner-finite} (proof of Theorem 1.6), it
   is quite elementary.
  However, I could not find it in the literature.
  Apparently, a natural generality here
  are spaces $H\setminus G/T$, where $G$ is a classical (or semisimple group),
  $H$ is a symmetric subgroup, and $T$ is the maximal unipotent subgroup in a parabolic
  subgroup.
  
  A possibility to describe double coset space implies a question about 
  harmonic analysis for $L^2$ on $\T_+(\alpha_1,\dots,\alpha_p)\setminus \GL(n)$.
  Such an analysis is possible, see Addendum to this paper (but it is not directly related 
  to the description of double cosets).

\section{Proof of Theorem \ref{th:}}

\COUNTERS 

{\bf \punct Proof of Proposition \ref{pr:}.a.}
We must show that $\chi_{ij}$ does not depend on the choice of a representative 
of a double coset. 

Reformulate the definition of $\chi_{ij}(A)$ in the following way. 
We consider the intersection
$$
H=H(A):=A^{-1} (W_j\oplus\dots \oplus W_q)\cap (V_1\oplus \dots \oplus V_i)
$$
and send 
$$
H(A)\to Z:= (V_1\oplus \dots \oplus V_i)\bigoplus (W_j\oplus\dots \oplus W_q)
$$
by the formula 
\begin{equation}
h\mapsto (h,Ah).
\label{eq:hAh}
\end{equation}
 Next, we send $Z$ to the quotient
$$
Z\Bigl/\Bigl(
(V_1\oplus \dots \oplus V_{i-1})\bigoplus (W_{j+1}\oplus\dots \oplus W_q)\Bigr)
\simeq V_i\oplus W_j.
$$
Clearly, the relation $\chi_{ij}(A)$ is the image of $H(A)$ in this space.

\sm 

Now, let $C\in \T_+$, consider $AC$ instead of $A$. Then $H(AC)=C^{-1} H(A)$.
The set of vectors $(h,Ah)$, see \eqref{eq:hAh} changes
to $(C^{-1}h,Ah)$. But $C^{-1}$ acts trivially in 
$$
V_1\oplus \dots \oplus V_{i}\Bigl/ V_1\oplus \dots \oplus V_{i-1}\simeq V_i.
$$
Therefore $\chi_{ij}(AC)=\chi_{ij}(A)$.

\sm 

{\bf \punct Proof of \eqref{eq:chiij}.}
To be definite, verify that $\ker \chi_{ij}=\dom \chi_{i(j+1)}$.
Write the equation
$$
\begin{pmatrix}
0\\ \vdots\\ 0\\ y_{j+1}\\\vdots \\ y_q
\end{pmatrix}
=A 
\begin{pmatrix}
x_1\\ \vdots\\ x_i\\0\\ \vdots \\0
\end{pmatrix}.
$$
A vector $x_i$ is contained in $\ker \chi_{ij}$
if there are $x_1$, \dots, $x_{i-1}$, $y_{j+1}$, \dots, $y_q$
such that  this equation is satisfied. This implies that
$(x_i, y_{j+1})\in \chi_{i(j+1)}$. In particular $x_i\in\dom \chi_{i(j+1)}$.

Conversely, let $x_i\in \dom \chi_{i(j+1)}$. Then there are
$x_i$ and
$x_1$, \dots, $x_{i-1}$, $y_{j+2}$, \dots, $y_q$ satisfying the equation.
This implies that $x_i\in \indef \chi_{ij}$.

\sm

{\bf\punct Verification  of \eqref{eq:kraj1}--\eqref{eq:kraj2}.}
To be definite, let us prove the statements from the first row (\ref{eq:kraj1}).

\sm 

Chow that $\indef\chi_{1j}=0$. Let $\eta\in\indef \chi_{1j}$. Then there $y_{j+1}$, \dots, $y_q$
such that
$$
\begin{pmatrix}
\vdots\\0\\ \eta\\ y_{j+1}\\ \vdots
\end{pmatrix}
= A \begin{pmatrix}
0\\0\\ \vdots
\end{pmatrix}.
$$
But the right hand side must be zero and $\eta=0$.

\sm 

The statement $\im \chi_{qj}=V_j$ follows from the surjectivity
of $A$.

\sm 

{\bf \punct The action of $\prod \GL[W_j]\times \prod\GL[V_i]$ on the double cosets space.}
Let $G$ be a group, $K$, $L$ its subgroups, $\wt K$ and $\wt L$ the normalizers
of $K$ and $L$. Then the group $\wt K/K\times \wt L/L$ acts
on $K\setminus G/L$.  Indeed, for $\kappa\in \wt K$, $\lambda\in \wt L$
we have
$$
\kappa^{-1}\cdot K \cdot g \cdot L\cdot \lambda=
\cdot K \cdot \kappa^{-1} g \lambda \cdot L,
$$
So this transformation sends double cosets to double cosets.
Clearly, orbits of $\wt K/K\times \wt L/L$ on $K\setminus G/L$
are in one-to-one correspondence with double cosets
$\wt K\setminus G/\wt L$.

In our case the normalizers of $\T_-(\dots)$ and $\T_+(\dots)$
are the groups $\P_-(\dots)$ and $\P_+(\dots)$, the quotients
are $\prod \GL[W_j]$ and $\prod \GL[V_i]$.

So let us describe the double coset spaces
\begin{equation}
\P_-[W_1,\dots,W_q]\setminus\GL(n)/ \P_+[V_1,\dots,V_p].
\label{eq:PGLP}
\end{equation}
Our subgroups contain the usual subgroups of  lower and upper
triangle matrices respectively. Applying the usual Gauss reduction
we observe that any double coset contains a 0-1-matrix%
\footnote{A 0-1 matrix is a matrix consisting of zeros and unit and containing
only one unit in each column and each row}, i.e., an element
of the symmetric group $S(n)$. After this reduction 
we can permute basis elements in each $V_i$ and in each $W_j$,
so the double coset space \eqref{eq:PGLP} is in one-to-one 
correspondence with
$$
\prod_{j=1}^q S(\beta_j) \setminus S(n)/\prod_{i=1}^p S(\alpha_i)
$$
Our matrix has a natural decomposition into $pq$ blocks, it
is important only the number of units in each block.
We formulate our observation in the following complicate form.

\begin{lemma}
	\label{l:0-1}
	For any double coset \eqref{eq:PGLP} there are canonical decomposition
	of each $V_j$ and each $W_i$ into a direct sum of coordinate subspaces
	\begin{align}
	V_i=\k^{\alpha_i}=\oplus_{j=1}^q \k^{\alpha_i^j}=:\oplus_{j=1}^q V_i^j;
	\\
	W_j=\k^{\beta_j}=\oplus_{i=1}^p \k^{\beta_j^i}=:\oplus_{i=1}^p W_j^i,
	\end{align}
	such that
	$$\alpha_i^j=\beta_j^i.$$
	A representative of the double coset is the map sending
	each $V_i^j$ to the corresponding $W_j^i$ coordinate-wise.
\end{lemma}

Less formally, we get a matrix of the form
\begin{equation}
J[\{V_i^j\},\{W_j^i\}]=
\left( 
\begin{array}{ccc|ccc|ccc|ccc}
\1&0&0 & 0&0&0 & 0&0&0 & 0&0&0\\
0&0&0 & \1&0&0 & 0&0&0 & 0&0&0\\
0&0&0 & 0&0&0 & \1&0&0 & 0&0&0\\
0&0&0 & 0&0&0 & 0&0&0 & \1&0&0\\
\hline
0&\1&0 & 0&0&0 & 0&0&0 & 0&0&0\\
0&0&0 & 0&\1&0 & 0&0&0 & 0&0&0\\
0&0&0 & 0&0&0 & 0&\1&0 & 0&0&0\\
0&0&0 & 0&0&0 & 0&0&0 & 0&\1&0\\
\hline 
0&0&\1 & 0&0&0 & 0&0&0 & 0&0&0\\
0&0&0 & 0&0&\1 & 0&0&0 & 0&0&0\\
0&0&0 & 0&0&0 & 0&0&\1 & 0&0&0\\
0&0&0 & 0&0&0 & 0&0&0 & 0&0&\1\\
\end{array}
\right).
\label{eq:BIG}
\end{equation}
Here $p=4$, $q=3$. We present a decomposition of a matrix
into blocks corresponding to the decompositions $\oplus V_i$ and
$\oplus W_j$, and refined blocks corresponding to decompositions
$\oplus \oplus V_i^j$ and $\oplus \oplus W_j^i$. Units are put in bold
to make them visible among zeros.

\begin{lemma}
	\label{l:hinge-canonical1}
For the matrix \eqref{eq:BIG}
the corresponding linear relations $\chi_{ij}$ are the following:
\begin{align}
\ker \chi_{ij}=V_i^{j+1}\oplus \dots\oplus V_i^p,\qquad
\dom \chi_{ij}=V_i^{j}\oplus \dots\oplus V_i^p;
\\
\indef \chi_{ij}=W_j^1\oplus \dots \oplus W_j^{i-1},\qquad
\im\chi_{ij}=W_j^1\oplus \dots \oplus W_j^{i}.
\end{align}
The operator 
$$\Theta(\chi_{ij}):\dom\chi_{ij}/\ker\chi_{ij}\to \im \chi_{ij}/\indef\chi_{ij}$$
is the identical map $V_i^j\to V_j^i$.	
\end{lemma}

We say that such a bi-hinge is {\it standard} and denote it by
$$
\hinge[\{V_i^j\},\{W_j^i\}]
$$

The statement is obvious.
Indeed, let us write $\chi_{32}$ for the matrix \eqref{eq:BIG} (the general case differs from considerations below only by  longer notation). We  apply this matrix to  a vector
$$
\left(
\begin{array}{ccc|ccc|ccc|ccc}
x_1^1&x_1^2&x_1^3& x_2^1&x_2^2&x_2^3& \xi^1&\xi^2&\xi^3&0&0&0
\end{array}
\right)
$$
and get
\begin{equation}
\left(
\begin{array}{cccc|cccc|cccc}
x_1^1&x_2^1&\xi^1&0 & x_1^2&x_2^2&\xi^2 & 0& x_1^3&x_2^3&\xi^3&0
\end{array}
\right).
\label{eq:long1}
\end{equation}
On the other hand this must be equal to
\begin{equation}
\left(
\begin{array}{cccc|cccc|cccc}
0&0&0&0 & \eta^1&\eta^2&\eta^3&\eta^4 & y_3^1& y_3^2& y_3^3& y_3^4
\end{array}
\right).
\label{eq:long2}
\end{equation}
Recall that the linear relation $\chi_{32}$ consists of vectors
$$\Bigl(\begin{pmatrix}
\xi^1&\xi^2&\xi^3
\end{pmatrix},
\begin{pmatrix}
\eta^1&\eta^2&\eta^3&\eta^4 
\end{pmatrix}
\Bigr)\in V_3\oplus W_2,
$$
for which there are $x_{\dots}^{\dots}$, $y^{\dots}_{\dots}$
such that \eqref{eq:long1} equals to \eqref{eq:long2}.
We get 
\begin{align}
\xi_1=0,\quad \xi_2=\eta_3,\quad \xi_3=y_3^3
\label{eq:var11}
\\
x_1^2=\eta_1,\quad x_2^2=\eta_2, \quad \eta_4=0.
\label{eq:var12}
\end{align}
The remaining equations contain no information:
\begin{align}
x_1=0,\quad x_2^1=0,\quad 0=0;
\label{eq:var21}
\\
x_1^3=y_3^1,\quad x_2^3=y_3^2, \quad 0=y_3^4.
\label{eq:var22}
\end{align}
The sets of variables in \eqref{eq:var11}--\eqref{eq:var12}
and \eqref{eq:var21}--\eqref{eq:var22} do not intersect (since we start with 0-1-matrix),
On the other hand the second system \eqref{eq:var21}--\eqref{eq:var22}
has a solution, since each variable is present in it only one time
(again, this is a priori clear, because we start with 0-1-matrix).

Finally, we get the linear relation $\chi_{23}$ consisting of vectors
$$
\Bigl(\begin{pmatrix}
	0&\xi^2&\xi^3
\end{pmatrix},
\begin{pmatrix}
	\eta^1&\eta^2&\eta^3&0
\end{pmatrix}
\Bigr),\qquad \text{where $\xi_2=\eta_3$.}
$$

{\bf \punct The action of $\prod\GL[V_i]\times \prod \GL[W_j]$ on the set of bi-hinges.}
Clearly, the group $\prod\GL[V_i]\times \prod \GL[W_j] $
acts on the space
$$V_1\times \dots \times V_p\times W_1\times \dots\times W_q,$$
therefore it acts on the set of bi-hinges.

\begin{lemma}
	\label{l:hinge-canonical2}
{\rm a)}	Any orbit of the group $\prod\GL[V_i]\times \prod \GL[W_j]$ on the set $$\Hinge[V_1,\dots,V_p;W_1,\dots, W_q]$$
contains a unique standard bi-hinge $\hinge[\{V_i^j\};\{W_j^i\}]$
{\rm (}as in
	Lemma {\rm \ref{l:hinge-canonical1})}.
	
	\sm 
	
{\rm b)} The stabilizer  $G[\{V_i^j\};\{W_j^i\}]\subset 
\prod_{i=1}^p\GL[V_i]\times  \prod_{j=1}^q \GL[W_j]$
 of a standard bi-hinge $\hinge[\{V_i^j\};\{W_j^i\}]$
is the semidirect 	of the reductive group
\begin{equation}
\prod_{i=1}^p\prod_{j=1}^q \GL[V_i^j]\simeq \prod_{j=1}^q\prod_{i=1}^p \GL[W^i_j]
\label{eq:reductive}
\end{equation}
and the unipotent group
\begin{equation}
\prod_{i=1}^p \T_-[V_i^1,\dots, V_i^q]
\times 
\prod_{j=1}^q \T_+[W_j^1,\dots,W_j^p].
\label{eq:unipotent}
\end{equation}
\end{lemma}

{\sc Proof.} a)
In a fixed $V_i$ we have a flag
$$
V_i=\dom\chi_{i1}\supset \ker \chi_{i1}=\dom\chi_{i2} \supset \ker \chi_{i2}=\dom\chi_{i3}\supset \dots \supset \ker\chi_{iq}=0. 
$$
We choose an element of $\GL[V_i]$ sending this flag to 
the flag of decreasing coordinate subspaces of the form
\begin{equation}
\begin{pmatrix}
0&\dots&0&*&\dots&*
\end{pmatrix}
,
\label{eq:**}
\end{equation} this flag is canonically determined by dimensions of 
$\ker \chi_{ij}$.

Similarly, we fix $W_j$, consider the flag
$$
0=\indef\chi_{1j}\subset\im \chi_{1j}=\indef\chi_{2j}\subset \im \chi_{2j}=
\indef \chi_{3j}\subset \dots\subset \im \chi_{pj}=W_j,
$$
and choose an element of $\GL[W_j]$ sending this flag to
a flag consisting of an increasing  sequence of subspaces
of the form
$\begin{pmatrix}
*&\dots&*&0&\dots&0
\end{pmatrix}$.

So for any element of our hinge we fixed positions
of its domain, kernel, image and indefinity.
After this fixing, it remains a possibility to choose coordinates
in subquotient of the flag \eqref{eq:**}.

We get the desired canonical form. The statement b) also becomes
obvious, since the stabilizer must regard the flags in each $V_i$ and $W_j$
and the maps $\Theta(\cdot)$.
\hfill $\square$   

\sm 

{\bf \punct Coincidence of stabilizers.} Thus we have a map 
$$
\T_-[V_1,\dots,V_p]\setminus \GL(n)/\T_+[W_1,\dots, W_q]\to \Hinge[V_1,\dots,V_p;W_1,\dots, W_q],
$$
which is
$\prod\GL[V_i]\times  \prod \GL[W_j]$-equivariant 
 and establishes a bijection of the sets
of orbits. For a proof of Theorem \ref{th:}
it is sufficient to show that this map establishes a bijective 
for each pair of corresponding orbits. So we must check the coincidence of stabilizers
of canonical representatives of orbits. So it suffices
to prove the following lemma

\begin{lemma}
Consider a representative of a double coset in the canonical
0-1-form	
$J[\{V_j^j\},\{W_j^i\}]$, see Lemma {\rm\ref{l:hinge-canonical1}},
and the corresponding standard hinge $\hinge[\{V_j^j\},\{W_j^i\}]$.
Then any element of $\prod \GL[V_i]\times \prod \GL[W_j]$
stabilizing the hinge stabilizes $J$.
\end{lemma}

{\sc Remark.} Notice that the inverse inclusion of stabilizers follows
from the equivariance. \hfill $\boxtimes$

\sm

{\sc Proof.} The stabilizer of a standard hinge is described
in Lemma \ref{l:hinge-canonical2}. It is a product of subgroups
\eqref{eq:reductive} and \eqref{eq:unipotent}. For the reductive factor
\eqref{eq:reductive} the statement is clear. The unipotent factor
 itself is a product, and it is sufficient  to prove the statement
 for any factor in \eqref{eq:unipotent}, say $\T_-[V_i^1,\dots,V_j^q]$.
 To be definite, we show what is happened for the matrix $J$ given by \eqref{eq:BIG}
 and a matrix 
 $$
 \begin{pmatrix}
 1&0&0\\
 x&1&0\\
 y&z&1
 \end{pmatrix}\in \T_+[V_2^1,V_2^2,V_2^3]\subset \GL[V_2].
 $$
 Multiplying \eqref{eq:BIG} by this element we get the matrix
$$
\left( 
\begin{array}{ccc|ccc|ccc|ccc}
\1&0&0 & 0&0&0 & 0&0&0 & 0&0&0\\
0&0&0 & \boxed{\1}&0&0 & 0&0&0 & 0&0&0\\
0&0&0 & 0&0&0 & \1&0&0 & 0&0&0\\
0&0&0 & 0&0&0 & 0&0&0 & \1&0&0\\
\hline
0&\1&0 & 0&0&0 & 0&0&0 & 0&0&0\\
0&0&0 & \mathbf{x}&\boxed{\1}&0 & 0&0&0 & 0&0&0\\
0&0&0 & 0&0&0 & 0&\1&0 & 0&0&0\\
0&0&0 & 0&0&0 & 0&0&0 & 0&\1&0\\
\hline 
0&0&\1 & 0&0&0 & 0&0&0 & 0&0&0\\
0&0&0 & \mathbf{y}&\mathbf{z}&\1 & 0&0&0 & 0&0&0\\
0&0&0 & 0&0&0 & 0&0&\1 & 0&0&0\\
0&0&0 & 0&0&0 & 0&0&0 & 0&0&\1\\
\end{array}
\right).
$$
We put nonzero symbols in bold to make them visible on the field of zeros.
Clearly, this matrix can be reduced to the initial form
\eqref{eq:BIG} by a left multiplication by an element
of $T_-[W_1,\dots, W_q]$. The boxed units allow to delete $x$, $y$, $z$.
\hfill $\square$

\section*{Addendum. The spaces $L^2$ on $T_+(\alphA)\setminus\GL(n,\C)$}

{\bf \boldmath A.1. The principal series of unitary representations
	 of the groups $\GL(n,\C)$.}
Denote by $\Lambda$ the set of pairs $\lambda|\lambda'$, where
$\lambda$, $\lambda'\in \C$, $\lambda-\lambda'\in \Z$, 
$\Im \lambda=\Im \lambda'$. In other words, we consider pairs
of the type
$$
\lambda|\lambda'=\tfrac{k+is}2\Bigl|\tfrac{-k+is}2, \qquad \text{where $k\in \Z$, $s\in \R$}.
$$ 
For $\lambda|\lambda'\in \Lambda$ we have a well defined 'generalized power'
of any nonzero $z\in \C$:
$$
z^{\lambda|\lambda'}=z^\lambda \ov z^{\,\lambda'}:=(z/|z|)^{k}\cdot |z|^{is}.
$$

Consider the subgroup $\B_+(n):=\P_+(1,\dots,1)\subset \GL(n,\C)$ consisting of all upper triangular
matrices,
\begin{equation}
C:=
\begin{pmatrix}
c_{11}&c_{12}&\dots&c_{1n}\\
0&c_{22}&\dots&c_{2n}\\
\vdots&\vdots&\ddots&\vdots\\
 0&0&\dots &c_{nn}
\end{pmatrix}
,
\tag{A.1}
\end{equation}
  A {\it signature} $\lambdA$ is a collection
 of the form
 $$\lambdA:=(\lambda_1|\lambda'_1,\dots,\lambda_n|\lambda'_n),
 \qquad \text{where $\lambda_j|\lambda'_j\in \Lambda$}.$$
 For such $\lambdA$ denote by $\chi_\lambdA(A)$ the character
 of $\B_+(n)$ defined by
 \begin{equation}
 \chi_\lambdA(A):=\prod c_{jj}^{\lambda_j|\lambda'_j}.
 \tag{A.2}
 \end{equation}
 By $\rho_\lambdA$ we denote the representation
 of $\GL(n,\C)$ unitary induced in the sense of Mackey (see, e.g., \cite{Kir}, Subsect.~13.2, \cite{BR}, Sect.~16, Sect.19.1) from a one-dimensional 
 representation $\chi_\lambdA$ of the subgroup $\B_+(n)$.
 Such representations are called representations of the {\it unitary
 	nondegenerate principal series}, see, e.g., \cite{BR}, Sect. 19.3.
 Representations $\rho_\lambdA$ and $\rho_\mU$
 are equivalent iff a collection $\{\mu_k|\mu_k'\}_{k=1,\dots,n}$ is obtained
 from a collection $\{\lambda_j|\lambda_j'\}_{j=1,\dots,n}$ by a permutation.
 Denote by $\Sigma_n$ the set of all signatures determined upto a permutation.
 
 Next, consider the action of $\GL(n,\C)\times\GL(n,\C)$ on $\GL(n,\C)$
 by left and right multiplications, $g\mapsto h_1^{-1}g h_2$. This determines 
 the left-right regular representation of $\GL(n,\C)\times\GL(n,\C)$ in 
 $L^2\bigl(\GL(n,\C)\bigr)$. According Gelfand and Naimark, see, e.g., \cite{BR},
 Sect. 14.4.A, this representation
 decomposes as a direct integral
 \begin{equation}
 \int_{\Sigma_n}\rho_\lambdA\otimes \rho_{-\lambdA} \,d\lambdA
 \tag{A.3}
 \end{equation}
 of representations of $\GL(n,\C)\times\GL(n,\C)$.
 
 \sm

 {\bf \boldmath A.2. The space $L^2(\T_+(\alpha_1,\dots,\alpha_p)\setminus\GL(n,\C)\bigr)$.}
 Let us decompose a quasiregular representation of $\GL(n,\C)$ in $L^2(\T_+(\alpha_1,\dots,\alpha_p)\setminus\GL(n,\C)\bigr)$.
 This problem has an additional symmetry.
 Indeed, the subgroup $\P_+(\alpha_1,\dots,\alpha_p)$ normalizes
  $\T_+(\alpha_1,\dots,\alpha_p)$. Therefore the quotient group
  $$
  \P_+(\alpha_1,\dots,\alpha_p)/\T_+(\alpha_1,\dots,\alpha_p)\simeq\prod_{i=1}^p \GL(\alpha_i,\C)
  $$
 acts on $\T_+(\alpha_1,\dots,\alpha_p)\setminus\GL(n,\C)$ by left multiplications.
 So we get a unitary representation of the group
 \begin{equation*}
 G:=
 \prod_{i=1}^p \GL(\alpha_i,\C) \times \GL(n,\C)
 \tag{A.4}
 \end{equation*}
in our $L^2$. 

For  elements $\lambdA^{j}\in \Sigma_{\alpha_j}$ denote by 
$$\lambdA^{1}\sqcup \dots\sqcup \lambdA^p\in \Sigma_n$$
the row obtained by concatenation of rows $\lambdA^{j}$.

\medskip

 {\bf {Theorem} A.1.}
{\it The decomposition of $L^2(\T_+(\alpha_1,\dots,\alpha_p)\setminus\GL(n,\C)\bigr)$
under the action of the group	 $G=\prod_{i=1}^p \GL(\alpha_i,\C) \times \GL(n,\C)$
is multiplicity-free and has the form}
$$
\int\limits_{\lambdA^1\in\Sigma_{\alpha_1}}\dots \int\limits_{\lambdA^p\in\Sigma_{\alpha_p}}
\Bigl(\rho_{-\lambdA^1}\otimes \dots \otimes \rho_{-\lambdA^p}\Bigr) \otimes \rho_{\lambdA^{1}\sqcup \dots\sqcup \lambdA^p}\, d\lambdA^p\dots d\lambdA^1.
$$

\medskip  

{\sc Proof.} We have  a space homogeneous with respect to the group $G$. The stabilizer
$G_0$
of the initial point consists of tuples
$$
\left\{
b_1\in \GL(\alpha_1), \dots, b_p\in \GL(\alpha_p), 
\begin{pmatrix}
b_1&a_{12}&\dots&a_{1n}\\
0&b_2&\dots&a_{2n}\\
\vdots&\vdots&\ddots&\vdots\\
0&0&\dots&b_{p}
\end{pmatrix}\in \P_+(\alpha_1,\dots,\alpha_p)
\right\}.
$$
 By definition our representation is induced from the trivial representation
 of the stabilizer $G_0$. Consider a lager group 
 $$G_0^*=\prod \GL(\alpha_i,\C)\times \P_+(\alpha_1,\dots,\alpha_p). 
 $$
 We apply inducing in stages (see, e.g., \cite{Kir}, Subsect.13.1,  \cite{BR}, Sect. 16.2),
 first from $G_0$ to $G_0^*$, second from $G_0^*$ to $G$.
 
 On the first step we have the same normal subgroup $\T_+(\alpha_1,\dots, \alpha_p)$ 
 in both $G_0$ to $G_0^*$, since the initial representation of $G_0$ is trivial, the induced representation
 is trivial on $\T_+(\alpha_1,\dots, \alpha_p)$. In fact we have the induction from
 $$
 G_0/
 \T_+(\alpha_1,\dots, \alpha_p)\qquad \text{to}\qquad G_0^*/\T_+(\alpha_1,\dots, \alpha_p).
 $$
 The second group is the double $\prod\GL(\alpha_j,\C)\times\prod\GL(\alpha_j,\C)$, the first
 group is $\prod\GL(\alpha_j,\C)$  embedded to the double as the diagonal. So the induced representation
 is the left-right representation of the double. Clearly, it is equivalent to the tensor product
 of the left-right regular representations of the factors $\GL(\alpha_j,\C)\times\GL(\alpha_j,\C)$,
 $$
 L^2\Bigl(\prod\GL(\alpha_j,\C)\Bigr)=\bigotimes_j L^2\bigl(\GL(\alpha_j,\C)\bigr)
 .
 $$
  We decompose spaces $L^2(\GL(\alpha_j,\C))$
 according (A.3).
 
In this way, we come to a direct integral of irreducible representations of $G_0^*$
having the form
\begin{equation}
\Bigl(\rho_{-\lambdA^1}(b_1)\otimes\dots\otimes  \rho_{-\lambdA^p}(b_p) \Bigr)\otimes \sigma_{\lambdA^1,\dots,\lambdA^p},
\tag{A.5}
\end{equation}
where $\sigma_{\lambdA^1,\dots,\lambdA^p}$ is the representation of $\P_+(\alpha_1,\dots,\alpha_p)$,
which is trivial on $\T_+(\alpha_1,\dots,\alpha_p)$ and is defined by
$$
\sigma_{\lambdA^1,\dots,\lambdA^p} \begin{pmatrix}
a_{11}&a_{12}&\dots&a_{1n}\\
0&a_{22}&\dots&a_{2n}\\
\vdots&\vdots&\ddots&\vdots\\
0&0&\dots&a_{pp}
\end{pmatrix}=\rho_{\lambdA^1}(a_{11})\otimes \dots \otimes \rho_{\lambdA^p}(a_{pp}).
$$
Notice, that $\sigma_{\lambdA^1,\dots,\lambdA^p}$ itself is an induced
representation. It is induced from a one-dimensional  representation of $\B_+(n)$, see
(A.1), namely from the character given by the formula
\begin{multline}
\zeta(C)=
\prod_{j=1}^{\alpha_1} c_{jj}^{\lambda_j^1\bigl|(\lambda_j^1)'}
\cdot 
\prod_{j=1}^{\alpha_2} c_{(\alpha_1+j)(\alpha_1+j)}^{\lambda_j^2\bigr|(\lambda_j^2)'}
\cdot
\prod_{j=1}^{\alpha_3} c_{(\alpha_1+\alpha_2+j)(\alpha_1+\alpha_2+j)}^{\lambda_j^3\bigr|(\lambda_j^3)'}
\cdot\, \dots=\\=
\chi_{\lambdA^{1}\sqcup \dots\sqcup \lambdA^p}(C),
\tag{A.6}
\end{multline}
where $\chi\dots(C)$ is determined by (A.1).

 Next, we consider the representation of $G$ induced from an irreducible representation
 (A.5)
 of $G_0^*$.  Since the factor $\prod \GL(\alpha_j)$ is present in both groups $G_0^*$ and $G$,
 actually we have the induction from $\P_+(\alpha_1,\dots,\alpha_p)$ to
 $\GL(n)$ (formally, we can refer to \cite{BR}, Sect.16.2.D, Theorem 3). But the representation $\sigma_{\lambdA^1,\dots,\lambdA^p}$ itself
 is induced from the character (A.6) of $\B_+(n)$. Applying inducing in stages
 we get that our representation of $\GL(n)$ is induced from the character (A.6)
 of $\B_+(n)$. But this a representation $\rho_{\lambdA^{1}\sqcup \dots\sqcup \lambdA^p}$ of the principal series.
 \hfill $\square$

\tt
\noindent
 Math. Dept., University of Vienna \\
\&Institute for Theoretical and Experimental Physics (Moscow); \\
\&MechMath Dept., Moscow State University;\\
\&Institute for Information Transmission Problems;\\
yurii.neretin@univie.ac.at
\\
URL: http://mat.univie.ac.at/$\sim$neretin/

\end{document}